\documentclass[12pt]{amsart}

\newcounter{defcounter}
\usepackage{amssymb,amsmath,amscd,graphicx, enumitem,
latexsym,amsthm}
\usepackage{amssymb,latexsym,eufrak,amsmath,amscd,graphicx}
  \usepackage[all]{xy}
  \usepackage{pdfsync}
\theoremstyle{plain}
\newtheorem{theorem}{Theorem}
\newtheorem{proposition}[theorem]{Proposition}

\newtheorem{lemma}[theorem]{Lemma}

\newtheorem{proposition.definition}[theorem]{Proposition/Definition}

\newtheorem{theoremalpha}{Theorem}

\newtheorem{conjecture}[theorem]{Conjecture}

\theoremstyle{definition}
\newtheorem{definition}[theorem]{Definition}

\newtheorem{remark}[theorem]{Remark}
\newtheorem{example}[theorem]{Example}

\newtheorem{question}[theorem]{Question}

\newtheorem{problem}[theorem]{Problem}

\newtheorem{definitionalpha}[theoremalpha]{Definition}

\newcommand{\lra}{\longrightarrow}
\newcommand{\lla}{\longleftarrow}
\newcommand{\noi}{\noindent}
\newcommand{\PP}{\mathbf{P}}

\newcommand{\CC}{\mathbf{C}}

\newcommand{\OO}{\mathcal{O}}

\newcommand{\II}{\mathcal{I}}

\newcommand{\frakm}{\mathfrak{m}}

\newcommand{\HH}[3]{H^{{#1}} \big( {#2} , {#3}
\big) }

\newcommand{\coker}{\textnormal{coker}}

\newcommand{\pr}{\prime}

\newcommand{\gon}{\textnormal{gon}}

\newcommand{\Tor}{\textnormal{Tor}}
\newcommand{\Sym}{\textnormal{Sym}}

\newcommand{\ol}[1]{\overline{#1}}

\numberwithin{theorem}{section}

\begin{document}

\title[Syzygies of  Projective Varieties of Large Degree]
{Syzygies of  Projective Varieties of Large Degree: \\ Recent Progress and Open Problems}

 \author{Lawrence Ein}
  \address{Department of Mathematics, University Illinois at Chicago, 851 South Morgan St., Chicago, IL  60607}
 \email{{\tt ein@uic.edu}}
 \thanks{Research of the first author partially supported by NSF grant DMS-1501085.}

 \author{Robert Lazarsfeld}
  \address{Department of Mathematics, Stony Brook University, Stony Brook, New York 11794}
 \email{{\tt robert.lazarsfeld@stonybrook.edu}}
 \thanks{Research of the second author partially supported by NSF grant DMS-1439285.}

\maketitle

 \section*{Introduction}
 
The purpose of this paper is to survey some recent work concerning the asymptotic behavior of the defining equations and higher syzygies of a smooth projective variety as the positivity of the embedding line bundle grows.

 To set the stage, we start with some rough history. Classically,  there was  interest in trying to say something about the equations defining suitably positive embeddings  of projective varieties. For example, let $C$ be a smooth projective curve of genus $g$, and let $L = L_d$ be a line bundle of degree $d \ge 2g+1$, giving rise to an embedding
\[   C \ \subseteq \ \PP H^0(L) \ = \ \PP^r, \]
where $r  = r_d = d-g$. Castelnuovo and others proved that $C$ is projectively normal, and cut out by quadrics as soon as $d \ge 2g + 2$.\footnote{There are actually three possible meanings for the statement that a projective variety $X \subseteq \PP^r$ is cut out by quadrics. The weakest is to ask that this simply be true set-theoretically. A more substantial condition is that $X$ be defined as a subscheme of $\PP^r$ by equations of degree two, ie that the twisted ideal sheaf $\II_{X/\PP^r}(2)$ be globally generated. The strongest possibility is that the homogeneous ideal $I_X$ of $X$ is generated by elements of degree two. All of the results described here hold in this last sense, although this isn't always the setting in which they were originally established.} 
   Mumford and his school studied the analogous (but much less elementary) questions for an abelian variety $A$ of arbitrary dimension. Specifically, consider an ample divisor $\Theta$ on $A$, and put $L = L_d = \OO_A(d\Theta)$. Then $L_d$ is very ample when $d \ge 3$, and  it defines a projectively normal embedding 
 in which $A$ is cut out by quadrics when $d \ge 4$. These issues were popularized in  \cite{Mumford1}, where Mumford also established that starting with any smooth projective variety $X$, a sufficiently positive Veronese re-embedding of $X$ is always cut out by quadrics.

In the early 1980s, as a byproduct of his work \cite{Kosz1} on Koszul cohomology, Mark Green realized that   results  of this type should be seen as the first cases of a much more general picture involving higher syzygies. Specifically, consider a very ample line bundle $L$ on a smooth projective variety $X$, defining an embedding
\[   X \, \subseteq \,  \PP H^0(L) \, = \, \PP^r, \]
where $r = r(L) = h^0(L) - 1$. 
Write $S = \Sym \, H^0(L)$ for the homogeneous coordinate ring of $\PP^r$, and put 
\begin{equation} 
\label{Graded.Ring.L}
R \ = \ R(X; L) \ = \ \oplus\,  H^0(X, mL). \notag\end{equation}
Thus $R$ is a finitely generated graded $S$-module, and so admits a minimal free resolution $E_\bullet = E_\bullet(X;L)$
\begin{equation}\label{Resolution}
 \xymatrix{
0 & R \ar[l]& E_0 \ar[l]  & E_1  \ar[l]  & \ar[l] \ldots & \ar[l] E_r \ar[l]  &  \ar[l]0 ,}
 \end{equation}
where
\[    E_p \ = \ \oplus S(-a_{p,j}). 
\]
Observe that $L$ is normally generated if and only if $E_0 = S$, in which case the remainder of $E_\bullet$ determines a minimal resolution of the homogeneous ideal $I_X \subseteq S$ of $X$. It is elementary that 
\[    a_{p,j} \, \ge \, p +1 \ \ \text{for all }\ j. \]
Green realized that the way to generalize the classical results is to ask when the first few terms of the resolution are generated in lowest possible degree.

The following definition formalizes Green's insight:
\begin{definitionalpha} \label{NpDef}
For $k \ge 0$ we say that $L$ satisfies Property ($N_k$) if $L$ defines a projectively normal embedding, and if
\[  E_p \ = \ \oplus \, S(-p-1) \ \text{ for } \ 1 \le p \le k. \ \qed\]
\end{definitionalpha}
\noi Thus ($N_0$) holds for $L$ if and only if $L$ is normally generated, and ($N_1$) is equivalent to requiring that in addition the homogeneous ideal $I_X$ of $X$ be generated by quadrics. The first non-classical condition is ($N_2$), which asks that if one chooses quadratic generators $Q_\alpha \in I_X$, then the module of syzygies among the $Q_\alpha$ should be spanned by relations of the form
\[  \sum \, L_\alpha \cdot Q_\alpha \ = \ 0, \]
where the $L_\alpha$ are \textit{linear} polynomials. For example, the resolution of the ideal of the rational normal cubic curve  $C \subseteq \PP^3$ has the shape
\[  0 \lla I_C \lla S(-2)^3 \lla S(-3)^2 \lla 0, \]
and so ($N_2$) holds. On the other hand, an elliptic quartic curve $E \subseteq \PP^3$ is a complete intersection of two quadrics, whose ideal is resolved by a Koszul complex:
\[  0 \lla I_E \lla S(-2)^2 \lla S(-4)\lla 0. \]
So in this case ($N_1$) holds but not ($N_2$). 

Green showed that the result of Castelnuovo et. ~al. ~on defining equations of curves admits a very natural generalization to higher syzygies:
\begin{theoremalpha}[Green, \cite{Kosz1}]\label{Green.Thm.Curves}
Let $L = L_d$ be a line bundle of degree $d$ on a smooth projective curve $C$ of genus $g$. If \[ d \, \ge \, 2g + 1 + k,\] then $L$ satisfies Property $(N_k)$. 
\end{theoremalpha}
This  result generated a great deal of interest and further work, much of it in the direction of finding extensions to other classes of varieties. For example, Green treated the case of Veronese embeddings in \cite{Kosz2}:
\begin{theoremalpha} \label{Veronese.Green}
The line bundle $\OO_{\PP^n}(d)$ satisfies $(N_k)$ for $d \ge k$. 
\end{theoremalpha} 
\noi Inspired by a conjecture of Mukai, Theorem \ref{Veronese.Green}  was generalized by the authors to arbitrary non-singular xvarieties in \cite{SAD}:
\begin{theoremalpha}
Let $X$ be a smooth projective variety of dimension $n$,  let $B$ and $P$ be respectively a very ample and a nef divisor on $X$. Then the line bundle \[  L_d \ = \ K_X + dB + P  \]  satisfies property $(N_k)$ provided that $d \ge n + 1 + k$. 
\end{theoremalpha}
\noi The case of toric varieties was studied in \cite{HSS}, and Galligo and Purnaprajna established  interesting results for surfaces in \cite{GP1}, \cite{GP2}, and \cite{GP3}. Arguably the deepest result along these lines is due to Pareschi \cite{Pareschi}, who extended the work of Mumford et.\,al. ~on abelian varieties to higher syzygies:
\begin{theoremalpha}[Pareschi, \cite{Pareschi}] Let $A$ be an abelian variety of arbitrary dimension $n$,  let $\Theta$ be an ample divisor on $A$, and put $L_d = \OO_A(d \Theta)$. If 
\[  d \ \ge \ k + 3, \]
then  $(N_k)$ holds for $L_d$.
\end{theoremalpha} 
\noi Pareschi's  argument used ideas involving the Fourier-Mukai transform, which were in turn systematized and extended in a very interesting series of papers by Pareschi and Popa \cite{PP1}, \cite{PP2}.  
Syzygies of abelian varieties were revisited from the viewpoint of local positivity in \cite{HwangTo}, \cite{LPP} and \cite{KL}. 

It is suggestive to summarize these results as asserting that Property ($N_k$) holds linearly in the positivity of the embedding line  bundle. More precisely, let $X$ be a smooth complex projective variety of dimension $n$, let $A$ and $P$ denote respectively an ample and an arbitrary divisor on $X$, and put
\begin{equation} L_d \ = \ dA + P. \end{equation}
Then one can recapitulate the results above by the following
\begin{theoremalpha} \label{Metatheorem.Statement}
There exist positive constants $C_1$ and $C_2$ depending on $X, A$ and $P$, such that $L_d$ satisfies property $(N_k)$ for 
\[ k \ \le \ C_1d  +C_2.\]
\end{theoremalpha}

This gives a good overall picture of the situation for curves of large degree, but when $\dim X = n \ge 2$ these results ignore most of the syzygies that can occur. Specifically, recall that the length of the resolution \eqref{Resolution} associated to a line bundle $L_d$ is essentially \[ r_d \ = \ r(L_d) \ = \ h^0(L_d) -1. \]
On the other hand, by Riemann-Roch
\[    r_d \ \sim \ \textnormal{(Constant)}\cdot d^n . \]
Hence when $n \ge 2$, the picture given by  Theorem \ref{Metatheorem.Statement}
  leaves open the possibility that the overall shape of the resolution of $L_d$ for $d \gg 0$ is quite different than what one might expect by extrapolating from Green's theorem on curves. In fact, the first indication that this is the case was a result of Ottaviani and Paoletti \cite{OP} asserting  that while ($N_k$) holds linearly for Veronese embeddings (Theorem \ref{Veronese.Green}), it also \textit{fails} linearly:
\begin{theoremalpha}[Ottaviani--Paoletti, \cite{OP}] \label{OP.Nonvan}
Property $(N_k)$ fails for $\OO_{\PP^n}(d)$ when $k \ge 3d - 2$.\footnote{They also conjecture -- and prove in the case of $\PP^2$ -- that ($N_p$) holds for $p < 3d-2$.}
\end{theoremalpha}
\noi The body of work surveyed in the present paper arose in an effort to understand systematically the asymptotic behavior of the syzygies   for very positive embeddings of higher-dimensional varieties.

In \S 1 we discuss and illustrate the main asymptotic non-vanishing theorem, and we state some conjectures that would complete the overall picture. In \S 2, we turn to the particularly interesting case of Veronese varieties, where following \cite{QPNV} we explain a very simple proof of the main cases of  non-vanishing. Section 3 centers on some results and conjectures concerning the asymptotics of Betti numbers. Finally, we return to curves in \S4, and explain the proof of the gonality conjecture from \cite{Gonality} and   discuss briefly the extension in \cite{WOS} of this result to higher dimensions. 

We deal throughout with projective varieties over the complex numbers, and we take the customary liberties of confusing divisors and line bundles. The reader may  refer to \cite{Eisenbud} for a presentation of the algebraic perspective on syzygies. Limitations of space and focus prevent us from discussing the very fundamental work of Voisin \cite{Voisin1}, \cite{Voisin2} on Green's conjecture on the syzygies of canonical curves, as well as its further developments e.g.      in \cite{AF}. We refer for example to Beauville's expos\'e \cite{Beauville} for an overview of the question and Voisin's results.  

During the course of the work reported here we have profited from discussions with many colleagues, including Marian Aprodu, David Eisenbud, Daniel Erman, Gabi Farkas, Mihai Fulger, Milena Hering, G. Ottaviani,   B. Purnaprajna,  Claudiu Raicu, Frank Schreyer, Jessica Sidman, David Stepleton, Bernd Sturmfels, Claire Voisin, David Yang, and Xin Zhou.

\numberwithin{equation}{section}

\section{Non-vanishing for asymptotic syzygies}

We start by fixing notation. Until further notice, $X$ is a smooth complex projective variety of dimension $n$, and we put 
\[ L_d \, = \, dA + P, \]
where $A$ is an ample and $P$ an arbitrary divisor. We always suppose that $d$ is sufficiently large so that $L_d$ is very ample, defining an embedding
\[  X \ \subseteq \ \PP H^0(X, L_d) \ = \ \PP^{r_d}, \]
where $r_d = h^0(L_d) - 1$. As in the Introduction, we denote by $S = \Sym \, H^0(L_d)$ the homogeneous coordinate ring of $\PP^{r_d}$. One can then form the minimal graded free resolution  associated to the ring $R(X;L_d)$ determined by $L_d$, but it will be useful to consider a slightly more general construction.

Specifically, fix a line bundle $B$ on $X$, and set
\[  R\, = \, R(X,B; L_d) \ = \ \oplus_m \, H^0(X, B + mL_d). \]
This is in the natural way a finitely generated graded $S$-module, and so has a miminal graded free resolution $E_\bullet = E_\bullet(X, B;L_d)$ as in equation \eqref{Resolution}. 

\begin{example} \label{O1.Example}
Consider the embedding 
\[   \PP^1 \, \subseteq \, \PP^3 \ \ , \ \ [s,t] \mapsto [s^3, s^2t, st^2, t^3] \]
of $\PP^1$ as the twisted cubic, which is cut out by the three quadrics
\[ 
Q_1 \, = \, XZ - Y^2 \ \ , \ \ Q_2 \, = \, XW - YZ \ \ , \ \ Q_3 \, = \, YW - Z^2. \]
Taking $B = \OO_{\PP^1}(1)$, the resulting module $R$ over $S = \CC[X,Y,Z,W]$ has two generators $e , f \in R_0$   corresponding to $s, t \in H^0(\PP^1, \OO_{\PP^1}(1))$. These satisfy the relations
\[    Y e - Xf  \, = \, 0 \ , \ Z  e - Y f\, = \, 0 \ , \ W  e - Zf \,= \, 0, \]
and we find the resolution
\[   
 \xymatrix{
0 & R   \ar[l]& S^2 \ar[l] & & & S^3(-1) \ar[lll]  _{\left(\begin{smallmatrix}{\ \,}Y & {\ \, }Z &{\ \, }	 W\\-X & -Y &-Z\end{smallmatrix}\right)}& & \ar[ll]_{\left(\begin{smallmatrix}{\ }Q_3 \\ - Q_2 \\ {\ }Q_1\end{smallmatrix}\right)}S(-3) &  \ar[l]0 .}   \ \ \qed
\]
\end{example}

We now come to the basic:
\begin{definition} \textbf{(Koszul cohomology groups).}
Define \[
K_{p,q} \big ( X, B; L _d\big) \ = \ \Big \{
\parbox{2.6in}{\begin{center} minimal generators of $E_p(X,B;L_d)$ of degree $p + q$\end{center}}
\Big  \}.  
\]
  Thus $K_{p,q}(X, B;L_d)$ is a finite-dimensional vector space, and
\[  E_p(X, B;L_d) \ = \  \underset{q}{\bigoplus}   \ K_{p,q}(X, B;L_d) \, \otimes _\CC \,  S(-p-q). \]
\end{definition}
 We refer to elements of $K_{p,q}$ as $p^{\text{th}}$ syzygies of weight $q$. When $B = \OO_X $ -- as in the Introduction -- we write simply $K_{p,q}(X; L_d)$. We recall at the end of this section that $K_{p,q}$ can be computed the cohomology of a bigraded Koszul-type complex.  
 
\begin{example}
In the situation of Example \ref{O1.Example}, one has
\[   K_{0,0} = \CC^2 \ \ , \ \ K_{1,0} = \CC^3 \ \  , \ \ K_{2,1} = \CC, \]
while $K_{p,q} = 0$ for all other $(p,q)$. \qed
\end{example}

\begin{example}
Assume that $B = \OO_X$. Then $L_d$ satisfies Property ($N_k$) if and only if 
\begin{align*} K_{0,q}(X; L_d) \ &= \ 0 \ \ \text{for } q \ne 0 \\ K_{p,q}(X; L_d) \ &= 0 \ \ \text{for } q \ne 1 \ , \ 1 \le p \le k .  \ \ \qed
\end{align*}
\end{example}

\begin{example} \textbf{(Betti diagrams).} It is often suggestive to display the dimensions of the various $K_{p,q}$ in tabular form, with  rows indexed by the weight $q$ and the columns corresponding to relevant values of $p$. For instance, the resolution computed in Example \ref{O1.Example} is summarized in the table:
\begin{center}
   \begin{tabular}{c |@{\hspace{10pt}} c@{\hspace{20pt}} c@{\hspace{20pt}} c}
 {} & 0 &1 & 2\\
  \hline  
0 &    2 & 3 & -- \\ 
 1&  -- & -- & 1 \\ 
  \end{tabular}
\end{center}
It is customary to use a dash  to indicate a zero entry. Note that the grading conventions are such that two adjacent entries on the same row correspond to a map in the resolution given by a matrix of linear forms. \qed
\end{example}

Fixing $B$, we now turn to the question of which of the groups $K_{p,q}(X, B; L_d)$ are non-vanishing for $d \gg 0$. This problem is framed by the following result, which shows that the situation is completely controlled when $q = 0$ or $q \ge n+1$.
\begin{proposition} \label{Easy.Kpq} For $d \gg 0$:
\begin{itemize}
\item[$($i$)$.] $K_{p,q}(X, B; L_d) \, = \, 0 \ \text{for }
q \ge n+2.$
\vskip 5pt
\item[$($ii$)$.] $ K_{p,0}(X,B; L_d) \, \ne \, 0 \ \Longleftrightarrow \ p \, \le \, r(B)$.
\vskip 5pt
\item[$($iii$)$.]
$K_{p, n+1}(X,B;L_d) \ne 0$ if and only if \[
r_d - n - r(K_X -B)  \ \le \  p  \ \le \ r_d - n.\]
\end{itemize}
\end{proposition}
\noi We refer to \cite[\S 5]{ASAV} for the proof. Statement (i) follows easily from considerations of Castelnuovo - Mumford regularity, while (ii) and (iii) are established by combining  arguments of Green \cite{Kosz1} and Ottaviani--Paoletti \cite{OP}.

\begin{example} \textbf{(Green's Theorem).} When $X$ is a curve and $B = \OO_X$, the  Proposition implies Green's Theorem \ref{Green.Thm.Curves}, at least for $d = \deg(L) \gg 0$.  In fact, it follows from (ii) that 
$K_{p,2}(X;L) = 0$ when
\[   (d - g) - g \, > \, p, \]
and since in any event all $K_{p,q} = 0$ for $q \ge 3$, this means that ($N_k$) must hold for $k \le (2g + 1) -p$. This is essentially the argument by which Green established the result in \cite{Kosz1}. \qed
\end{example}

The main  non-vanishing theorem from \cite{ASAV} asserts that from an asymptotic perspective, essentially all of the remaining Koszul groups are non-zero.

\begin{theorem} \label{Asympt.Non.Van}
Fix $1 \le q \le n$. There exist constants $C_1, C_2 > 0$ $($depending on $X, B, A$ and $P$$)$ with the property that for $d \gg 0$,
\[   K_{p,q}(X, B; L_d) \, \ne \, 0 \]
for every value of $p$ with 
\begin{equation} \label{Range.of.p}  C_1 \cdot d^{q-1} \ \le p \ \le \ r_d \, - \, C_2 \cdot d^{n-1}. \end{equation}
\end{theorem}
\noi Some  effective statements appear in \cite{Zhou1} and in Theorem \ref{CM.Syzygies}   below. 

To get a feeling for the statement, fix $q \in [1, n]$ and set 
\[   w_q(d) \ = \ \frac{ \# \big \{ p \in [1, r_d] \mid K_{p,q}(X, B; L_d) \ne 0 \, \big \}}{\# \big \{ p \in [1, r_d]\big \} }, \]
so that $w_q(d)$ measures the proportion of potentially non-zero weight $q$ syzygies that are actually non-zero. Recalling that $r_d = O(d^n)$, the Theorem implies that 
\[  \lim_{d \to \infty} \, w_q(d) \ = \ 1. \]
In terms of the corresponding Betti diagram, one can visualize this as asserting that except for some negligibly small regions, the rows recording syzygies of weights $q = 1, \ldots, n$ are entirely filled by non-zero
entries.

The proof of the Theorem in \cite{ASAV} involves a rather complicated induction on dimension, the idea being that one can use suitable secant planes to produce non-zero syzygies. In the next section we will explain a much quicker argument for the case $X = \PP^n$ (or more generally when $X \subseteq \PP^N$ is projectively Cohen-Macaulay). However we would like to propose a heuristic explanation, which however we've never been able to  push through.

Taking $B = \OO_X$ for simplicity,  fix a hypersurface $\ol{X} \subset X$  and consider the embedding of $\ol{X}$ defined by $\ol{L}_d$ defined by the restriction of $L_d$. This gives rise to a commutative diagram:
\begin{equation} \label{Divisor.Diagram} \begin{gathered}
\xymatrix @R=.25pc  { X & \subseteq & \PP^{r_d} \\
\rotatebox{90}{$\subseteq$} & & \rotatebox{90}{$\subseteq$}\\
\ol{X} & \subseteq & \PP^{{\ol{r}_d}}}
\end{gathered}
\end{equation}
where
\[   r_d \, = \, r(X, L_d) \, = \, O(d^n) \ \ , \ \ \ol{r}_d \, = \, r(\ol{X},  \ol{L}_d) \, = \, O(d^{n-1}). \]
Now we can consider $\ol{X}$ as a subvariety both of $\PP^{\ol{r}_d}$ and $\PP^{r_d}$, and it is elementary that (roughly speaking):
\small
\[  
\big( \text{Syzygies of } \ol{X}\subset \PP^{{r}_d}\big ) \  = \ \big( \text{Syzygies of } \ol{X}\subset \PP^{\ol{r}_d}\big ) \otimes \big( \text{Koszul resolution of } \PP^{\ol{r}_d} \subseteq \PP^{r_d} \big) 
\]
\normalsize
(see \cite[\S1]{Kosz2} for the precise statement). By induction on dimension one can suppose that syzygies of many different weights appear in the resolution of $\ol{X}$ in $\PP^{\ol{r}_d}$,  and then the same will be true of the resolution of $\ol{X}$ in $\PP^{r_d}$ thanks to the presence of the very large Koszul complex appearing on the right. One \textit{expects} that this should finally force many non-vanishing Koszul groups in the resolution of $X$, but unfortunately it is not clear to us how to rule out the  (unlikely) vanishing of various maps in the long exact sequence relating the syzygies of $X$ and of $\ol{X}$. 

\begin{remark} \textbf{(Stanley-Reisner ideals of subdivisions).} The Stanley-Reisner ideal $I_{\Delta} $ of a simplicial complex $\Delta$ is a monomial ideal in a polynomial ring that encodes the combinatorics of $\Delta$. In their interesting paper \cite{CKW}, Conca, Juhnke-Kubitzke and Welker study the asymptotics of the syzygies of the ideals associated to repeated subdivisions of a given complex $\Delta$. They find that these satisfy the same sort of picture as occurs in the geometric setting: almost all of the Betti numbers that could be non-zero are in fact non-zero. \qed
\end{remark}

Returning to the situation of Theorem \ref{Asympt.Non.Van}, it is natural to ask what happens for those values of $p$ outside the range governed by the statement. We conjecture that the lower bound appearing in \eqref{Range.of.p} is actually the best possible in the sense that one has vanishing of $p^\text{th}$ syzygies for smaller $p$.

\begin{conjecture}  \textnormal{\textbf{(Asymptotic vanishing)}}. \label{NonVanConj}
Fix $q \in [2, n]$. In the situation of Theorem \ref{Asympt.Non.Van}, there is a constant $C_3$ $($depending on $X$, $A, B$ and $P)$ such that
\[   K_{p,q}(X, B:L_d) \ = \ 0 \ \text{ for } p \, \le \, C_3 \cdot d^{q-1}\]
when $d\gg 0$. 
\end{conjecture}
\noi When $q = 2$, this essentially follows from Theorem \ref{Metatheorem.Statement}  (which remains valid in the presence of an arbitrary twisting divisor $B$). In some situations -- for example for the Veronese embeddings discussed in the next section -- one can verify that the Conjecture is valid when $q = n$ (see Example \ref{Top.Weight.Rmk}). In general, Raicu \cite{Raicu} shows that knowing the conjecture for $X = \PP^n$ implies its truth for arbitrary varieties. We consider the Conjecture to be the main open problem concerning the rough asymptotics of the $K_{p,q}$.

Finally, we recall how in practice one   computes the $K_{p,q}(X, B; L_d)$. Writing $L$ in place of $L_d$, the basic result is the following:
\begin{proposition} \label{Koszul.Cohom.Prop} The group
$K_{p,q}(X, B; L)$ is the cohomology of the Koszul-type complex
\Small
\begin{gather*} {} \\ 
 \rightarrow  \Lambda^{p+1}H^0(L) \otimes  H^0\big((q-1)L + B\big)\rightarrow  \Lambda^{p}H^0(L) \otimes  H^0\big(qL + B\big)  \rightarrow\Lambda^{p-1}H^0(L) \otimes  H^0\big((q+1)L + B\big) \rightarrow .\\
\end{gather*}
\normalsize
Here the differential
\[  \Lambda^{p}H^0(L) \otimes  H^0\big(qL + B\big)  \lra \Lambda^{p-1}H^0(L) \otimes H^0\big( (q+1)L +B \big) \]
is given by 
\[ ( s_1 \wedge \ldots \wedge s_p) \otimes t \,\mapsto \, \sum (-1)^i (s_1 \wedge \ldots \wedge \widehat{ s_i} \wedge \ldots \wedge s_p) \otimes s_i \cdot t.\] 
\end{proposition}
\noi This is essentially a reflection of the symmetry of $\Tor$. In brief, write $\underline \CC = S/S_+$ for the quotient of the polynomial ring $S$ by the irrelevant maximal ideal. Tensoring a miminal resolution of $R = R(X, B; L)$   by $\underline{\CC}$, one sees that
\[   K_{p,q}(X, B; L) \ = \ \Tor_p(R, \underline{\CC})_{p+q}.\]
On the other hand, one can also compute these Tor's starting from the Koszul resolution of $\underline{\CC}$ and tensoring by $R$, and this leads to the complex appearing in the Proposition.

\section{Veronese Varieties}
In this section we discuss the particularly interesting case of Veronese varieties, where one can obtain effective statements. Specifically, we aim to establish non-vanishings for the groups
\[  
K_{p,q}(n,b;d) \ =_\text{def} \ K_{p,q}\big(\PP^n, \OO_{\PP^n}(b);\OO_{\PP^n(d)} \big)\]
for fixed $b$ and large $d$. Note that
\[
K_{p,q}(n,b;d) \ = \  K_{p,q+1}(n,b-d;d),
\]
so there is no harm in assuming that $0\leq b \leq d-1$.

The main result here, which was established in \cite{ASAV}  (for a slightly smaller range of the parameters) and much more quickly in \cite{QPNV}, is the following:
\begin{theorem}  \label{VeroneseNonVanishing}
Fix  $b\ge 0$ and $q \in [0,n]$. Then 
\[
K_{p,q}(n,b;d) \ \ne \ 0 \]
for any   \[ d \ \ge \ b + q + 1\] and 
all $p$ in the range
\small
\[
\binom{q+d}{q}- \binom{d-b-1}{q}-q\  \leq \ p\ \leq \ \binom{n+d}{n} - \binom{d+n-q}{n-q} + \binom{n+b}{q+b}-q-1. \tag{*}
\]
\normalsize 
\end{theorem}
\noi When $b = 0$ this result was established independently by Weyman.
\vskip 10pt
\begin{example}
Take $n = 2$ and $b = 0$. Then the Theorem asserts that  $K_{p,2}(\PP^2; \OO_{\PP^2}(d)) \ne 0$ for
\[
3d - 2 \, \le \, p \, \le \, \binom{d+2}{2} - 3, 
\]
which is exactly the result of Ottaviani--Paoletti \cite{OP} cited in the Introduction (Theorem \ref{OP.Nonvan}). 
\end{example}

We believe that the Theorem actually gives the precise non-vanishing range for Veronese syzygies:
\begin{conjecture} \label{Veronese.Conjecture}
In the situation of the Theorem, one has 
\[   K_{p,q}(n,b;d) \ = \ 0 \]
when $p$ lies outside the range $(*)$.
\end{conjecture}
\noi Our belief in the conjecture stems in part from the fact that several quite different approaches to the non-vanishing lead to exactly the same numbers, as well as from the following:

 \begin{remark} \textbf{(Syzygies of extremal weights).} \label{Top.Weight.Rmk}
 At least for $d \gg 0$, it is established in \cite[Remark 6.5]{ASAV} that the result is best-possible for $q = 0$ and $q = n$, i.e. that $K_{p,0}(n, b; d) = 0$ and $K_{p,n}(n,b;d) = 0$ when $p$ lies outside the stated range. This provides at least some  evidence for Conjectures \ref{Veronese.Conjecture}
 and \ref{NonVanConj}. 
 \end{remark}

 Theorems \ref{Asympt.Non.Van} and \ref{VeroneseNonVanishing} suggest that the syzygies of a given variety become 	quite complicated as the positivity of the embedding grows. In the case of Veronese varieties, one can try to make this more precise via representation theory. Specifically, 
the groups $K_{p,q}(n;d)$ are representations of $\textnormal{SL}(n+1, \CC)$, and hence decompose into irreducibe representations. It is then natural to ask about the shape of this decomposition  as $d \to \infty$. For fixed $p$  the $K_{p,q}$ vanish for large $d$ when $q \ge 2$ by virtue of  Theorem \ref{Veronese.Green}, so the  interesting case is that of $K_{p,1}$. One expects the  precise decomposition of $K_{p,1}$ to run up against essentially intractable questions of plethesym, but one can hope to get a picture of the complexity of these groups by counting the number of irreps that appear. This is given by a very nice result of Fulger and Zhou:
 \begin{theorem}[Fulger--Zhou, \cite{FZ}]
 Fix $p$, and an integer $n \ge p$. Then as $d \to \infty$ the groups $K_{p,1}(n;d)$ contain exactly on the oder of $d^p$ irreducible representations  of $\textnormal{SL}(n+1, \CC)$.
 \end{theorem}
\noi They also show that one gets on the order of $d^{(p^2+p)/2}$ irreps counting  multiplicities.  Thus  Veronese syzygies do indeed become quite complicated from a representation-theoretic perspective. The idea of Fulger and Zhou is to constuct a convex polytope whose lattice points parametrize the highest weights of representations appearing in  these Koszul groups.  We note that some related results appear in \cite{Rubei}.

\begin{remark} \textbf{(Toric varieties).} In his  interesting paper \cite{Zhou2}, Zhou studies  the distribution of torus weights for all the $K_{p,q}$ on a toric variety $X$. For a given range of $p$, Zhou describes explicitly the closure of the weights (suitably normalized), which may or may not fill out the polytope defining the toric projective embedding of $X$. \qed\end{remark}

Theorem \ref{VeroneseNonVanishing} was established in \cite{ASAV} by keeping track  of the secant constructions used in that paper, but a much quicker proof appears in \cite{QPNV} which reduces the question to some elementary computations with monomials. 
Write $S = \CC[z_0, \ldots, z_n]$ for the homogeneous coordinate ring of $\PP^n$. In view of  Proposition \ref{Koszul.Cohom.Prop}, the groups $K_{p,q}(n,b; d)$ that we are interested in are the cohomology of the complex
\[
... \lra \Lambda^{p+1} S_d \otimes S_{(q-1)d+b} \lra  \Lambda^{p} S_d \otimes S_{qd + b} \lra  \Lambda^{p-1} S_d \otimes S_{(q+1)d + b} \lra ...
\]
In principle one could hope to prove the non-vanishing of these groups by simply   writing down explicitly a suitable cocycle, but we do not know how to do this. 

However consider  the ring
\[  \ol{S} \ = \ S / (z_0^d, \ldots, z_n^d) \ . \]
We think of $\ol{S}$ as the algebra spanned by monomials in which no variable appears with exponent $\ge d$, with multiplication governed by the vanishing of the $d^{\text{th}}$ power of each variable. Now since $z_0^d, \ldots, z_n^d$ forms a regular sequence in $S$, the dimensions of the Koszul cohomolgy groups of $\ol{S}$ are the same as those of $S$, ie $K_{p,q}(n,b;d)$ is isomorphic to the cohomology of the complex
\[
... \lra \Lambda^{p+1} \ol{S}_d \otimes \ol{S}_{(q-1)d+b} \lra  \Lambda^{p} \ol{S}_d \otimes \ol{S}_{qd + b} \lra  \Lambda^{p-1} \ol{S}_d \otimes \ol{S}_{(q+1)d + b} \lra ...  \ . 
\]
Here the presence of many zero-divisors enables one easily to exhibit  non-vanishing cohomolgy classes.

We illustrate how this works by proving the Ottaviani--Paoletti statement (Theorem \ref{OP.Nonvan}) that \[K_{3d-2,2}(2;d) \ \ne \ 0 \] provided that $d \ge 3$. 
Writing (to lighten notation) $\ol{S} = \CC[x,y,z]/(x^d, y^d,z^d)$, this is equivalent by what we have just said to showing that the complex 
\[
 \Lambda^{3d-1} \ol S_d \otimes \ol S_{d} \lra  \Lambda^{3d-2} \ol S_d \otimes \ol S_{2d} \lra  \Lambda^{3d-3} \ol S_d \otimes \ol S_{3d}    \  . \tag{*}
\]
has non trivial homology. To this end, note first that if $m_1, \ldots, m_{3d-2}$ are any monomials of degree $d$ that are each divisible by $x$ or $y$, then the element
\[  c  \, = \,   m_1 \wedge \ldots \wedge m_{3d-2} \otimes x^{d-1}y^{d-1}z^2 \ \in \  \Lambda^{3d-2} \ol S_d \otimes \ol S_{2d} \tag{**} \]
is a cycle for (*). It remains to show that by chosing the monomials $m_i$ suitably we can arrange that $c$ is not a boundary. We will achieve this by taking the $m_i$ to be all the factors of $x^{d-1}y^{d-1}z^2$.

Specifically,  observe  that $x^{d-1}y^{d-1}z^2$ has exactly $3d-2$ monomial divisors of degree $d$ with exponents $\le d-1$, viz:
\begin{gather*}
x^{d-1}y \, , \, x^{d-2}y^2\, , \, \ldots\, , \,, x^2y^{d-2}\, , \,x y^{d-1}     \\
x^{d-1}z\, , \, x^{d-2}y z \, , \, \ldots\, , \, xy^{d-2}z\, , \, y^{d-1}z   \\
x^{d-2}z^2 \, , \, x^{d-3}yz^2 \, , \, \ldots \, , \, xy^{d-3}z^2 \, , \, y^{d-2}z^2.
\end{gather*}
We claim that if we use these as the $m_i$ in (**), then the resulting cycle $c$ represents a non-zero cohomology  class. In fact, suppose that $c$  were to appear even as a term in the Koszul boundary of an element
\[  e \ = \ n_0  \wedge   n_1   \ \ldots \  \wedge n_{3d-2} \otimes g,\]
where the $n_i$ and $g$ are monomials of degree $d$. After re-indexing   we can suppose that 
\[   c \ = \ n_1 \wedge \ldots \wedge n_{3d-2} \otimes n_0g. \] Then the $\{ n_j \}$   with $j \ge 1$ must be a re-ordering of the monomials $\{ m_i \}$ dividing $x^{d-1}y^{d-1}z^2$. On the other hand $n_0g = x^{d-1}y^{d-1}z^2$, so $n_0$ is also such a divisor. Therefore $n_0$ coincides with one of $n_1, \ldots, n_{3d-2}$,  and hence $e  = 0$, a contradiction.  Observe that if $m_{3d-1}, \ldots, m_p$ are additional  monomials that annihilate $x^{d-1}y^{d-1}z^2$ in $\ol{S}$, then the same argument shows that
\begin{equation} \label{cocycle} \big( m_1 \wedge \ldots \wedge m_{3d-2} \wedge m_{3d-1} \wedge \ldots \wedge m_p \big) \otimes x^{d-1}y^{d-1}z^2\end{equation} represents a non-zero class in $K_{p,2}(2;d)$, and in fact different choices of $m_{3d-1}, \ldots, m_p$ yield linearly independent classes.

With more careful book-keeping, it turns out that this approach gives exactly the statement appearing in Theorem \ref{VeroneseNonVanishing}. In fact, a similar argument yields an effective statement  analogous to Theorem \ref{VeroneseNonVanishing} for the Koszul cohomology groups of any projectively Cohen-Macaulay variety $X \subseteq \PP^N$ of dimension $n$:
\begin{theorem} \label{CM.Syzygies}
Denote by $c(X)$ the Castelnouvo-Mumford regularity of $\OO_X$, and put
\[
r_d \,= \,h^0(X, \OO_X(d)) \ \ , \ \ r^\pr_d \,= \, r_d - (\deg X)(n+1).
\]
Then for $q \in [1, n-1]$, and  $d \ge b + q + c(X) + 1$:
\[ K_{p,q}(X,\OO_X(b); \OO_X(d))\  \ne \ 0\] for every value of $p$ satisfying
\small
\[
\deg(X) (q + b + 1) \binom{d + q-1}{q-1} \ \le \ p \ \le \  {r}^\pr_d - \deg(X) (d-q-b) \binom{d + n-q-1 }{n-q-1}.
\]
\normalsize
\end{theorem}
\noi Analogous statements hold, with slightly different numbers, when $q =0$ and $q = n$

\section{Betti numbers}

In this section we discuss some results and conjectures from \cite{ARBT} concerning the asymptotics of  the Betti numbers of a very positive embedding. We keep notation as above: so $X$ is a smooth projective variety of dimension $n$, and we consider for large $d$ the embedding $X \subseteq \PP^{r_d}$ defined by the complete linear series associated to the line bundle
\[ L_d \ = \ dA + P. \]
Given a twisting line bundle $B$, and   weight $q \in [1,n]$, we will be interested in the  dimensions\[   k_{p,q}(X, B; L_d) \ =_{\text{def}}\ \dim  K_{p,q}(X,B; L_d) \ \ , \ \ k_{p,q}(X, B; L_d) \ =_{\text{def}}\ \dim K_{p,q}(X; L_d) \]
as functions of $p$ for $d \gg 0$.

\begin{figure}  
\includegraphics[scale=.8]{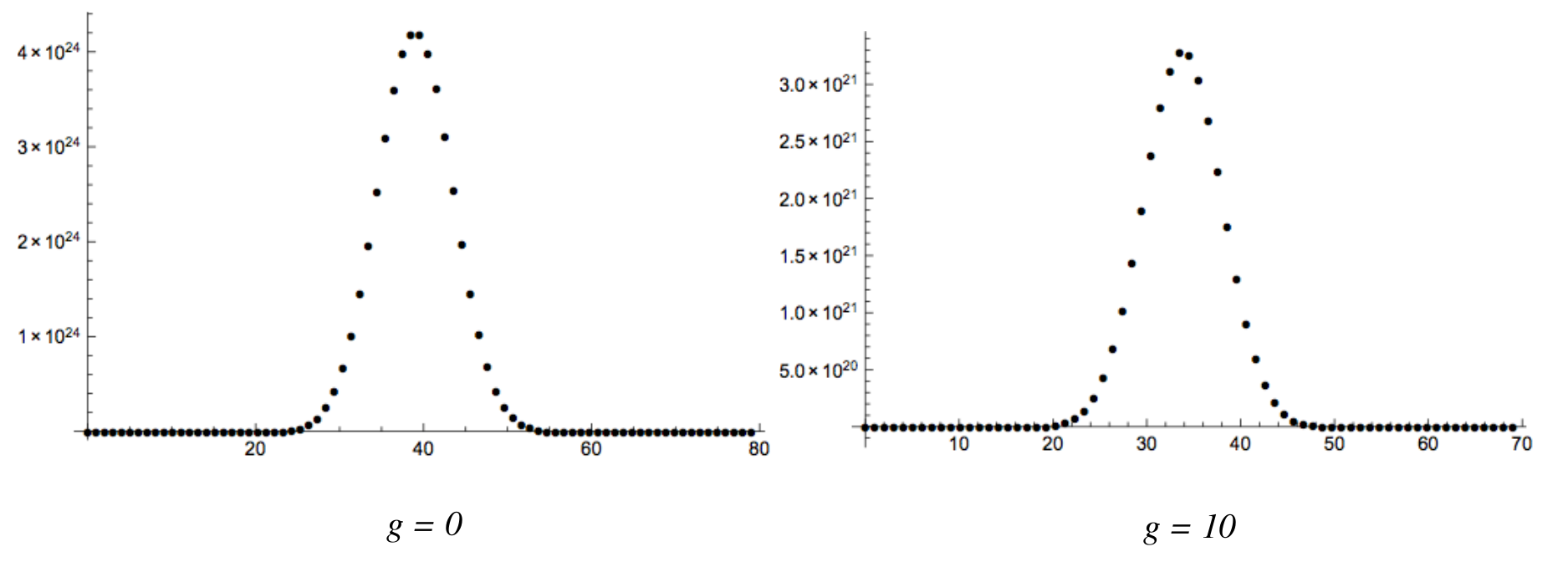}
\caption{Plots of $k_{p,1} $ for bundles of degree $80$ on curves of genus $g = 0$ and $ 10$}
\label{PlotCurve}
\end{figure}

The first case to consider is that of curves. Here Theorem \ref{Green.Thm.Curves} implies that for all except $g$ values of the parameter $p$, only weight one syzygies occur. In these instances  $k_{p,1}$ can be computed as an Euler characteristic, and one finds that for $p \le r_d - g = d - 2g$:
\[
k_{p,1}(X; L_d)  \ = \ \binom{r_d}{p} \left( \frac{-pd}{r_d} + (r_d+1) - \frac{d+1-g}{p+1} \right) .
 \]
The dominant term here is the binomial coefficient: Figure \ref{PlotCurve} shows plots of the $k_{p,1}(X;L_d)$ for a line bundle  of degree $ d = 80$ on  curves  of genus $0$ and $10$.\footnote{As we shall see in the next section, on a curve $X$ of genus $g$  the last $g$  Betti numbers $k_{p,1}$ for $d-2g \le p \le r_d = d-g$ depend on the intrinsic geometry of $X$ when $g \ge 3$. However the variation is small compared to the value of $k_{p,1}$ for $p \approx \frac{r_d}{2}$, and so is not visible graphically.}
More precisely, it follows from Stirling's formula that the function $k_{p,1}(X; L_d)$ becomes Gaussian as $d \to \infty$ in the following sense:
\begin{proposition} \label{Betti.Number.Curves}
Choose a sequence $\{ p_d \}$  of integers such that 
\[ p_d \ \to\ \frac{r_d}{2} + a \cdot \frac{\sqrt{r_d}}{2}\]
for some fixed number $a$ \textnormal{(}ie.\ $\lim_{d \to \infty} \,  \frac{2p_d - r_d}{\sqrt{r_d}} = a$\textnormal{)}.  Then as $d \to \infty$,
\[   
\left ( \frac{1}{2^{r_d} }\cdot \sqrt{\frac{2\pi}{r_d}}\, \right ) \cdot k_{{p_d},1}(X;L_d)  \to e^{-a^2 / 2}. \]
\end{proposition}

We conjectue that the same pattern holds universally:
\begin{conjecture} \label{Normal.Betti.Conjecture}
Returning to a smooth projective variety $X$ of dimension $n$, fix $q \in [1,n]$. Then there is a normalizing function $F_q(d)$ $($depending on $X$ and geometric data$)$ such that\[
F_q(d)\cdot k_{p_d,q}(X;L_d) \lra e^{-a^2 / 2} 
\]
as $d \to \infty$ and $p_d \to  \frac{r_d}{2} \ + \ a \cdot \frac{\sqrt {r_d}}{2} $.
\end{conjecture}
\noi One expects slightly more generally that the analogous statement is true for the dimensions $k_{p,q}(X, B; L_d)$ with $B$ a fixed twisting line bundle. 

It is not hard to establish  lower and upper bounds for the quantites in question that are  Gaussian in shape. For example, using the cocycles \eqref{cocycle} together with Proposition \ref{Koszul.Cohom.Prop} one sees that if $3d -2  < p < \binom{d+2}{2} -2$  then 
\[
\binom{\frac{d(d-1)}{2}}{p} \ \le \  k_{p,2}\big( \PP^2; \OO_{\PP^2}(d) \big ) \ \le \ \binom{\frac{(d+2)(d+1)}{2}}{p}\cdot (2d+1)(d+1),\]
but unfortunately the two bounds don't match up. In fact, the conjecture has not been verified for any single variety of dimension $n \ge 2$. 
 One could imagine that the large Koszul complex govering the embedding of $\PP^{\ol{r}_d} \subseteq \PP^{r_d}$ appearing in \eqref{Divisor.Diagram} comes into play here, but we don't have much of a picture how to make this precise. It would already be very interesting to have a conceptual -- rather than simply computational -- explanation for Proposition \ref{Betti.Number.Curves}. 

While the actual evidence in favor of Conjecture \ref{Normal.Betti.Conjecture} may seem skimpy, the main content of \cite{ARBT} was to   argue that at least the picture suggested by the conjecture is probabilistically very natural. To explain this in an especially simple setting, consider the Betti numners associated to   $B = \OO_{\PP^2}(-1)$ and $L_d = \OO_{\PP^2}(d)$ on $\PP^2$. In this case
\[    k_{p,q}(\PP^2,  B; L_d) \ =  \ 0  \ \text{ for } q \ne 1, 2, \]
ie the corresponding Betti table has only two rows. 
By the Boij-S\"oderberg theory of Eisenbud and Schreyer \cite{ES}, these Betti numbers can be expressed as non-negative linear combinations of those of certain ``pure modules." Specificlly, there exist modules $\Pi_i$ $(1 \le i \le r_d)$    having the property  
\[   K_{p,1}(\Pi_i) \ne 0 \, \Leftrightarrow \,  0 \le p < i \ \ , \ \ K_{p,2}(\Pi_i) \ne 0 \ \Leftrightarrow\  i  \le p \le r_d, \] 
together with rational numbers $x_i = x_i(\PP^2, B: L_d) \ge 0$ such that
\[   k_{p,q}(\PP^2, B; L_d) \ = \ \sum_{i = 0}^{r_d} \, x_i \cdot k_{p,q}(\Pi_i) \tag{*} \]
for all $p, q$.\footnote{Graphically, the Betti table  of $\Pi_i$ consists of $i$ non-zero entries in the $q = 1$ row, followed by $r_d - i$ non-zero entries in the $q=2$ row, with zeroes elsewhere.} We may call the $x_i$ the Boij-Soderberg coefficients of the Betti table of $B$ with respect to $L_d$. 

Now for arbitrary $ x_i \ge 0$, the right hand side of (*) defines the Betti numbers of a module with the given Boij-S\"oderberg coefficients,  which one might view as the potential Betti table of a surface. In order to test whether the behavior predicted by the conjecture is ``typical" or not, we ask what happens if we choose the $x_i$ randomly. By scaling one may suppose that $x_i \in [0,1]$, so consider  the hyper-cube $\Omega_r \ = \ [0,1]^r $ parametrizing $r$-tuples of Boij-S\"oderberg coefficients.  Given
\[
x \ = \ \{ x_i \} \ \in \ \Omega_r, \]
denote by 
\begin{equation} \label{InVitroKpq}
k_{p,q}(x) \ = \ \sum_{i = 0}^{r} \, x_i \cdot k_{p,q}(\Pi_i)  \end{equation}
the entries of the corresponding $2 \times r$ Betti table. 
Stated rather informally, we show that with high probability, the behavior predicted by the Conjecture holds for such a random Betti table:\begin{theorem} Fix $q = 1$ or $q =2$. 
Then as $r \to \infty$, with probability $ = 1$ the Betti numbers $k_{p,q}(x)$ satisfy the analogue of Conjecture \ref{Normal.Betti.Conjecture} when $x \in \Omega_r$ is sampled uniformly at random.
\end{theorem} \noi There is a similar statement for the random Betti tables modeling the syzygies of smooth varieties of dimensions $n \ge 3$. We refer to \cite{QPNV} for precise statements. It is also shown there that the statement is quite robust in the sense that the same conclusion holds if $x = \{ x_i \}$ is sampled with respect to many other probability measures on $\Omega_r$. 

\begin{remark} \textbf{(Asymptotic Boij-S\"oderberg coefficients).}
Returning to the global situation, we ask the following
\begin{quote} \textsc{Question.} Can one normalize the Boij-S\"oderberg coefficients $x_i(\PP^2, B; L_d)$ and the relevant values of $i$ in such a way that as $d \to \infty$ they arise as the values of a smooth function defined on a dense set in its domain? 
\end{quote}
Experience with asymptotic invariants of linear series suggests that something along these lines might well to be the case. The difference
$k_{p,1}(B; L_d) - k_{p-1,2}(B; L_d)$ can be computed as the Euler characteristic of a  vector bundle on $\PP^2$, and numerical experiments show that one gets good visual agreement with this difference if one takes the $x_i$ in \eqref{InVitroKpq}
  to be themselves  the values of a suitable Gaussian function. Of course one would like to ask the same question also in dimensions $n \ge 3$, but here there is some ambiguity in chosing Boij-S\"oderberg data. \qed
\end{remark}

\section{Asymptotic $K_{p,1}$ and the gonality conjecture}

The picture that we have discussed so far focuses on the rough overall structure of asymptotic syzygies, with statements largely independent of specific geometric hypotheses. However as observed \cite[\S 5]{ASAV}, one can hope for more precise results for the groups $K_{p,1}(X, B;L_d)$: in particular, for $d \gg 0$ one can expect that the values of $p$ for which these groups vanish to depend only on the geometry of $B$.  Results along these lines were established in   \cite{Gonality} and \cite{WOS}. The case of curves, treated in \cite{Gonality}, is particularly interesting as it leads to the proof of an old conjecture from \cite{GL}, so we start with this.  

Suppose then that $C$ is a smooth projective curve of genus $g$, fix a divisor $B$ on $C$, and  let $L_d$ be line bundle of degree $d \gg 0$ on $C$, so that $r_d = d -g$. Proposition \ref{Easy.Kpq} implies that if $d\gg 0$ then:
\begin{align*}
K_{p,0}(C, B; L_d) \, \ne \, 0\ \ &\Longleftrightarrow \  p\, \le \,r(B) \\
K_{p,2}(C, B; L_d) \, \ne \, 0\ \ &\Longleftrightarrow \ \ r_d - 1 - r(K_C -B) \, \le \, p \, \le \,r_d -1.
\end{align*} 
It follows that $K_{p,1}(C, B; L_d) \ne 0$ for
\[   r(B) + 1 \, \le \, p \, \le r_d - 2 - r(K_C -B), \]
since in this range none of the other $K_{p,q}$ appear. However this leaves open the 
\begin{question} \label{Kp1.Question}
\noi For which values of $p$ is $K_{p,1}(C, B; L_d) \ne 0$ when $d \gg 0$?
\end{question}
\noi Moreover by Serre duality \cite[\S 2c]{Kosz1} the groups
\begin{equation} \label{Dual.Groups} K_{p,1}(C, B; L_d) \ \ \text{and} \ \ K_{r_d - 1-p, 1}(C, K_C - B; L_d) \end{equation}
are dual, so it is enough to answer Question \ref{Kp1.Question} for $p \le r(B)$. 

The case $B =K_C$ was considered by Green in \cite{Kosz1}. These Koszul cohomology groups  control the syzygies of the so-called \textit{Arbarello--Sernesi module}
\[     M(C, K_C)  \ = \ \oplus \, H^0(C, K_C + mL_d). \]
When $g \ge 1$ this module has $g$ generators 
in degree $0$ corresponding to a basis of $H^0(C, K_C)$, and Green showed that $K_{0,1}(C, L_d) = 0$, i.e. that these generate $M$ as a module.\footnote{This is equivalent to the assertion that the multiplication map
\[  H^0(K_C) \otimes H^0(mL_d) \lra H^0(K_C + mL_d) \]
is surjective when $m \ge 1$, which is clear since $K_C$ is globally generated.}  More interestingly, he also showed that $K_{1,1}(C, K_C; L_d) = 0 $ -- in other words that the Arbarello-Sernesi module has a linear presentation -- if and only if $C$ is not hyperelliptic. It is natural to ask about the higher $K_{p,1}(C, K_C; L_d)$ for large $d$. 

Recall that the \textit{gonality} $\gon(C)$ of $C$ is by definition the least degree of a branched covering $C \lra \PP^1$. It is not hard to see that if 
$ \gon(C) \le  p+1$, then \[  K_{p,1}(C, K_C; L_d) \,\ne \, 0\] for large $d$.\footnote{A simple argument proceeds by noting that by duality, $K_{p,1}(C, K_C; L_d) \ne 0 $ if and only if $K_{r_d -1-p,1}(C; L_d) \ne 0$. But if $C$ carries a $g^1_{p+1}$, then for $d \gg 0$ the linear series in question sweeps out a rational normal scroll $\Sigma \subseteq \PP^{r_d}$ of dimension $p +1$ containing $C$, and the Eagon-Northcott resolution of the ideal of $\Sigma$ gives rise to the required weight one syzygies of $C$. }  
Motivated in part by his celebrated conjecture on the syzygies of canonical curves, this led Green and the second author (somewhat half-heartedly) to propose in \cite{GL} the
\begin{conjecture} \label{Gonality.Conjecture}
For $d \gg 0$, 
\[  K_{p,1}(C, K_C; L_d) \, \ne \, 0 \ \Longleftrightarrow \   \gon(C) \, \le \, p+1.\]
\end{conjecture}
\noi Drawing on Voisin's spectacular proof \cite{Voisin1}, \cite{Voisin2} of Green's conjecture for general canonical curves, Aprodu and Voisin \cite{Aprodu}, \cite{AproduVoisin} proved the Conjecture for many classes of curves, in particular for general curves of every gonality.

Recall that a line bundle $B$ on $C$ is said to be \textit{$p$-very ample} if for every effective divisor $ \xi \subseteq C$ of degree $p + 1$, the restriction map 
\[    H^0(C, B) \lra H^0(C, B \otimes \OO_{\xi}) \]
is surjective. Thus $B$ is $0$-very ample if and only if it is globally generated, and $B$ is $1$-very ample if and only if it is very ample. 
It follows from Riemann--Roch that the canonical bundle
 $K_C$ fails to be $p$-very ample if and only if 
 \[   
  \gon(C) \, \le \, p+1.\]
Therefore Conjecture \ref{Gonality.Conjecture} is a consequence of \begin{theorem} \label{Gon.Thm}
Fix a line bundle $B$ on $C$. Then 
\[  K_{p,1}(C, B; L_d) \, = \, 0 \ \text{ for } d \gg 0 \]
if and only if $B$ is $p$-very ample. 
\end{theorem}

\begin{remark} \textbf{(Resolution of curve of large degree).} It follows from the Theorem that one can read off the gonality of a curve $C$ from the resolution of the ideal of $C$ in any one embedding of sufficiently large degree. In fact, as in \eqref{Dual.Groups} the group $K_{p,1}(C, K_C; L_d)$ is dual to $K_{r_d -1 - p}(C; L_d)$. Therefore the gonality of $C$ is characterized as the least integer $c$ such that
\[  K_{r_d - c ,1}(C; L_d) \, \ne \, 0  \tag{*}\]
for any line bundle of degree $d \gg 0$. A result of Rathmann described in the next Remark shows that in fact it suffices here that $d \ge 4g - 3$.   Together with Proposition \ref{Easy.Kpq}, (*) means that one has a complete understanding of the grading of the resolution of the ideal of a curve of large degree.\qed
\end{remark}

\begin{remark}\textbf{(Rathmann's theorem).}
Rathmann \cite{Rathmann} has established an effective statement that essentially completes the story for curves. Specifically, he proves the following very nice
\vskip -10pt
\hskip 20pt\parbox[c]{5in}{
\begin{theorem} \label{Rathmann.Thm} Assume that $B$ is $p$-very ample, and that $L$ is any line bundle satisfying the vanishings
\[   H^1(C, L) \, = \, H^1(C, L-B) \, = \, 0. \]
Then $K_{p,1}(C, B; L) = 0$. 
\end{theorem}
}

\noi Thus for example Conjecture \ref{Gonality.Conjecture} holds for any line bundle $L_d$ of degree $d \ge 4g-3$. \qed
\end{remark}

 Theorem \ref{Gon.Thm} is surprisingly  quick and effortless to prove: like Poe's purloined letter, it turns out essentially to have been sitting in plain sight. The idea is to use Voisin's Hilbert schematic interpretation of syzygies, and reduce the matter to a simple application of Serre vanishing. Specifically, denote by $C_{p+1}$ the $(p+1)^\text{st}$ symmetric product of $C$, which we view as parameterizing effective divisors of degree $p+1$ on $C$. A line bundle $B$ on $C$ determines a vector bundle \[ E_B \ = \ E_{p+1, B} \] of rank $p+1$ on $C_{p+1}$, whose fibre at $\xi \in C_{p+1}$ is the $(p+1)$-dimensional vector space $H^0(C , B \otimes \OO_{\xi})$.  There is a natural evaluation map of vector bundles
 \begin{equation} \label{eval.map}  \textnormal{ev}_{B} \, : \, H^0(C, B) \otimes_{\CC} \OO_{C_{p+1}} \lra E_{p+1, B} \end{equation} which induces an isomorphism
 \[ \HH{0}{C_{p+1}}{E_B} \ = \ \HH{0}{C}{B}. \]
 Note that $\text{ev}_B$ is surjective as a map of bundles if and only if $B$ is $p$-very ample.
 
Given a line bundle $L$ on $C$, consider next the line bundle 
\[ \mathcal{N}_L \ = \ \mathcal{N}_{p+1,L} \ =_{\text{def}} \ \det E_{p+1, L} \]
on $C_{p+1}$. One can show that taking  exterior powers in the evaluation map \eqref{eval.map}  for $E_L$ gives rise to an isomorphism
\[    \HH{0}{C_{p+1}}{\mathcal{N}_L} \ = \ \Lambda^{p+1} H^0(C, L). \]
We now return to   \eqref{eval.map} and twist through by $\mathcal{N}_L$: using the computations of $H^0$ just stated, this gives rise to a homomorphism
\[    \HH{0}{C}{B} \otimes \Lambda^{p+1} \HH{0}{C}{L} \lra \HH{0}{C_{p+1}}{E_B \otimes \mathcal{N}_L }, \tag{*} \] 
and Voisin shows in effect that
\[  \HH{0}{C_{p+1}}{E_B \otimes \mathcal{N}_L } \ = \ Z_{p}(C,B;L) \]
is the space of cycles in the Koszul complex from Proposition \ref{Koszul.Cohom.Prop} computing $K_{p,1}(C,B;L)$.\footnote{Voisin actially worked on the universal family over the Hilbert scheme $C_{p+1}$, which is perhaps how this argument escaped notice.} Therefore $K_{p,1}(C,B;L) = 0$ if and only if the mapping (*) is surjective. 

Now asssume $B$ is $p$-very ample. Then $\text{ev}_B$ is surjective as a map of sheaves, and writing 
\[   M_B \ = \ M_{p+1, B} \ =_{\text{def}} \ \ker ( \text{ev}_B), \]
the vanishing of $K_{p,1}(C,B; L)$ will follow if we show that
\[  \HH{1}{C_{p+1}}{M_B \otimes \mathcal{N}_L} \ = \ 0  
\]
for $\deg (L) \gg 0$. 
But this is a consequence of
\begin{lemma} \label{Serre.Van.Lemma}
The line bundles $\mathcal{N}_L$ on $C_{p+1}$ satisfy Serre vanishing. More precisely, given any coherent sheaf $\mathcal{F}$ on $C_{p+1}$ there exists an integer $d_0 = d_0(\mathcal{F})$ with the property that
\[ \HH{i}{C_{p+1}}{\mathcal{F} \otimes \mathcal{N}_L} \ = \ 0 \ \ \text{ for all } \ i > 0
\] 
provided that $\deg(L) \ge d_0$. \qed
\end{lemma}
\noi To establish the more precise Theorem \ref{Rathmann.Thm}, Rathmann essentially replaces this appeal to Serre vanishing with a proof by descending induction on $q$ of an effective vanishing theorem for twists of $\Lambda^q M_B$.  The authors had used Griffiths vanishing to give a much weaker effecitve
 statement in \cite{Gonality}.
 
 \begin{remark} \textbf{(Growth of $\mathbf{k_{p,1}(C, B; L_d)}$).}
 The same setup yields some information about the dimension of $K_{p,1}(C,B;L_d)$ when $B$ is not $p$-very ample. In fact, put
 \[  \gamma_p(B) \ = \ \dim \, \big \{ \xi \in C_{p+1} \mid H^0(B) \lra H^0(B \otimes \OO_\xi) \ \text{ is not surjective} \big \}.\]
Applying Lemma \ref{Serre.Van.Lemma} to $\coker (\text{ev}_B)$ shows that if $d \gg 0$, then $k_{p,1}(B, C; L_d)$ is a polynomial in $d$ of degree $\gamma_p(B)$.  In his very interesting paper \cite{Yang}, Yang proves that on a smooth projective variety $X$ of arbitrary dimension, $\dim K_{p,1}(X, B; L_d)$ is a polynomial in $d$ for $d \gg 0$. \qed
 \end{remark}

\begin{remark} \textbf{(The secant conjecture).} The paper \cite{GL} proposed another conjecture that would interpolate between Green's Theorem \ref{Green.Thm.Curves} and his conjecture on canonical curves. Specifically, it was proposed that if $L$ is a $p$-very ample line bundle on a curve $C$ with
\[  \deg(L) \, \ge \, 2g + p + 1 - 2h^1(C,L) - \textnormal{Cliff}(L), \]
then $L$ satisfies Property ($N_p$). In their very nice paper \cite{FK}, Farkas and Kemeny prove this when $C$ and $L$ are general.  Kemeny carries this further in \cite{K}. \qed
\end{remark}

 It is natural to ask whether and in what form Theorem \ref{Gon.Thm} extends to higher dimensions. When $\dim X = n \ge 2$ there are two divergent notions of positivity for a line bundle $B$: $p$-very amplitude, which asks that $H^0(X, B) \lra H^0(X, B \otimes \OO_{\xi})$ be surjective for all subschemes of length $p+1$, and $p$-jet amplitude:
 \begin{definition}
 A line bundle $B$ on a smooth projective variety $X$ is said to be $p$-jet very ample of for every effective zero-cycle
 \[   w \, = \, a_1 x_1 + \ldots + a_s x_s \]
  of degree $p+1 = \sum a_i$ on $X$, tthe natural map
\[   \HH{0}{X}{B}  \lra \HH{0}{X}{B \otimes  \OO_X/ \frakm_w}\]
 is surjective, where 
 \[ \frakm_w \, =_{\text{def}} \, \frakm_1^{a_1} \cdot \ldots \cdot \frakm_s^{a_s}, \]
$\frakm_i  \subseteq \OO_X$ being the ideal sheaf of $x_i$. \end{definition}
\noi When $\dim X \ge 2$, this is a stronger condition than $p$-very amplitude.

Inspired by Yang's interpretation of Koszul cohomology in \cite{Yang}, Yang and the authors establish in \cite{WOS} the following:
\begin{theorem} Let $X$ be a smooth projective variety, and let $B$ be a line bundle on $X$. 
If  $B$ is   $p$-jet  very ample, then 
\[  K_{p,1}(X, B; L_d) \ = \ 0 \ \text{ for } \ d \gg 0. \]
Conversely, if there is a \textnormal{reduced} zero cycle $w = x_1 + \ldots + x_{p+1}$ that fails to impose independent conditions on $H^0(X,B)$, then
\[ K_{p,1}(X, B; L_d) \ \ne \ 0 \ \ \text{ for  all } \ d \gg 0. \]
\end{theorem}
\noi The first statement is proved by working on a cartesian self-product of $X$, establishing a vanishing of a group that contains the indicated $K_{p,1}$ as a summand.\footnote{This is an idea that goes back to Green in \cite{Kosz2}.}

To complete this picture, there remains:
\begin{problem}
Find necessary and sufficient conditions for the vanishing of $K_{p,1}
(X, B; L_d)$ when $d \gg 0$.
\end{problem}
\noi It does not seem out of the question that the failure of $B$ to be $p$-jet very ample in general implies the non-vanishing of this group. 
We had originally imagined the $p$-very amplitude of $B$ would control the matter, but a heuristic argument due to Yang casts some doubt on this possibility.\footnote{The arguments of Voisin in \cite{Voisin1} show that $K_{p,1}$ is computed by cohomolgy on the principal component of $\textnormal{Hilb}^{p+1}(X)$ parameterizing smoothable schemes, but when $\dim X \ge 3$ the failure of $B$ to be $p$-very amplitude  could be witnessed by a point on a different component of the Hilbert scheme. (Of course it's conceivable that the right condition involves smoothable schemes.)}

 %
 %
 %
 %

 \end{document}